\documentclass[microtype]{gtpart}     

\usepackage[all]{xy}

\title{Moduli space of nonnegatively curved metrics on manifolds of dimension $4k+1$}

\author{Anand Dessai}
\givenname{Anand}
\surname{Dessai}
\address{Department of Mathematics, University of Fribourg, Switzerland}
\email{anand.dessai@unifr.ch}
\urladdr{https://homeweb.unifr.ch/dessaia/pub/}


\arxivreference{https://arxiv.org//abs/2005.04741}

\def\wrt{with respect to}

\def\proofend{\hbox to 1em{\hss}\hfill $\blacksquare $\bigskip }

\newtheorem*{theorem*}{Main Theorem}

\newtheorem{theorem}{Theorem}[section]
\newtheorem{proposition}[theorem]{Proposition}
\newtheorem{lemma}[theorem]{Lemma}
\newtheorem{remark}[theorem]{Remark}
\newtheorem{remarks}[theorem]{Remarks}

\theoremstyle{definition}
\newtheorem{defn}[theorem]{Definition}

\def\Z{{\mathbb Z}}
\def\R{{\mathbb R}}
\def\Q{{\mathbb Q}}
\def\C{{\mathbb C}}
\def\N{{\mathbb N}}

\def\ccc {{\textbf{c}}}
\def\spinc{Spin ^\ccc}

\newcommand\sign{\mathrm{sign}}

\begin{document}

\begin{abstract}
In each dimension $4k+1\geq 9$ we exhibit infinite families of closed manifolds with fundamental group $\Z _2$ for which the moduli space of metrics of nonnegative sectional curvature has infinitely many path components. Examples of closed manifolds with finite fundamental group with this property were known before only in dimension $5$ and dimensions $4k+3\geq 7$.
\end{abstract}

\maketitle


\section{Introduction}\label{Section: Introduction}
We give examples of closed manifolds of dimension $4k+1$ with $k\geq 2$ for which the moduli spaces of metrics of nonnegative sectional curvature and positive Ricci curvature have infinitely many path components.

For a closed manifold $M$, let $\mathcal{R}_{sec\geq 0}(M)$ denote the space of Riemannian metrics of nonnegative sectional curvature on $M$ endowed with the smooth topology. The diffeomorphism group $\text{Diff}(M)$ acts on $\mathcal{R}_{sec\geq 0}(M)$ by pulling back metrics. The orbit space $
\mathcal{M}_{sec\geq 0}(M):=\mathcal{R}_{sec\geq 0}(M)/\text{Diff}(M)$
equipped with the quotient topology is called the {\em moduli space} of metrics of nonnegative sectional curvature on $M$. The corresponding notation will be used for the moduli space of metrics satisfying other curvature bounds.

A basic problem in Riemannian geometry is to determine whether a given manifold admits a metric with prescribed curvature properties. If this is the case, one may ask whether the respective moduli space carries some interesting topology.
In contrast to scalar curvature, where surgery techniques are available, little is known about the topology of moduli spaces of metrics satisfying lower bounds (nonnegative or positive) on sectional or Ricci curvature.

The first results in this direction are due to Kreck and Stolz \cite{KS93}, who introduced an invariant for certain $(4k+3)$-dimensional spin manifolds which is constant on path components of the moduli space of metrics of positive scalar curvature. Kreck and Stolz used this invariant to show that there exists an Aloff-Wallach space for which the moduli space of metrics of positive sectional curvature is disconnected. They also exhibited an infinite family of $7$-dimensional Witten manifolds for which the moduli spaces of metrics of positive Ricci curvature have infinitely many path components \cite{KS93}. Another more basic invariant to distinguish path components is the relative index of Gromov and Lawson \cite[p. 327]{GL83}. Using these invariants manifolds in dimension $4k+3\geq 7$ have been found (see \cite{KPT05,DKT18,D17,G17}) for which the moduli space $\mathcal{M}_{sec\geq 0}$ has infinitely many path components (see also \cite{W11,TW17,G20} as well as \cite{TW15} for related results).

Dessai and Gonz\'alez-\'Alvaro \cite{DG21} used relative $\eta$-invariants to show that for every homotopy $\R P^5$ the moduli space $\mathcal{M}_{sec\geq 0}$ has infinitely many path components (see also \cite{W20}).
Here we apply $\eta$-invariants to prove that manifolds with this property also exist in all dimensions $4k+1$ with $k\geq 2$.

\begin{theorem*} In each dimension $4k+1\geq 9$ there are infinitely many closed manifolds $M_i$, $i\in \N $, with pairwise nonisomorphic integral cohomology for which the moduli space $\mathcal{M}_{sec\geq 0}(M_i)$ of metrics of nonnegative sectional curvature has infinitely many path components. The same holds true for the moduli space $\mathcal{M}_{Ric > 0}(M_i)$ of metrics of positive Ricci curvature on $M_i$.  
\end{theorem*}
It follows that the corresponding spaces of metrics, $\mathcal{R}_{sec\geq 0}(M_i)$ and\linebreak $\mathcal{R}_{Ric > 0}(M_i)$, also have infinitely many path components.
 
In combination with
 \cite[Proposition 2.8]{BKS11} the theorem implies that for every such manifold the moduli space of complete metrics of nonnegative sectional curvature on the total space of a real line bundle over $M_i$ has infinitely many path components.

The manifolds in the theorem above may be described as total spaces of two-stage iterated fiber bundles over $\C P^1$ with fibers $\C P^{2k-1}$ and $S^1$ (see the next section for definitions and details) and are closely related to the manifolds considered in \cite{Wi81,WZ90,KS93,KPT05,DKT18}. They can also be described as quotients of the product of round spheres $S^3\times S^{4k-1}$ by a free isometric action of $S^1\times \Z _2$. The metrics which represent distinct path components in the respective moduli space are obtained as submersion metrics and have nonnegative sectional and positive Ricci curvature. To distinguish components, we compute relative $\eta$-invariants for these metrics. The construction can also be carried out for $k=1$, in which case one obtains a finite number of $5$-dimensional $\Z _2$-quotients of $S^2\times S^3$. Their moduli spaces of metrics of nonnegative sectional curvature and positive Ricci curvature also have infinitely many path components (see \cite{W20, GW22} and Remark \ref{dim five remark} for related results).

This paper is structured as follows.  In the next section we introduce a family of  $(4k+1)$-dimensional manifolds with fundamental group $\Z _2$ which are total spaces of two-stage iterated fiber bundles and collect some of their topological properties. In Section \ref{section diffeo finite}, we give a rough diffeomorphism classification for these manifolds. More precisely, we first study their homotopy type via Postnikov towers and then apply the exact surgery sequence to show that certain infinite subfamilies belong to only finitely many oriented diffeomorphism types.

The manifolds come with a submersion metric of nonnegative sectional and positive Ricci curvature, which, when lifted to the universal cover, extends in a nice way to an associated disk bundle. They also carry a $Spin^\ccc$-structure and a flat line bundle for which the relative $\eta $-invariant of the corresponding Dirac operator is nontrivial. This is explained in Sections \ref{section nonnegative curvature} and \ref{section spinc structures}. Computations for the relative $\eta$-invariant via equivariant index theory are detailed in Section \ref{section eta invariants}. These computations are then used in the final section to prove the main theorem.

\subparagraph*{Acknowledgements} I would like to thank Ian Hambleton for stimulating discussions and the referees for their valuable comments. Also many thanks to David Gonz\'alez-\'Alvaro for a useful remark on an earlier version of this paper. This work was supported in part by the SNSF-Project 200021E-172469 and the DFG-Priority programme SPP 2026.

\section{A family of $\mathbf{(4k+1)}$-dimensional manifolds}\label{section family}
In this section we describe a family of simply connected manifolds of nonnegative sectional curvature which will be used in the proof of the main theorem. A manifold in this family is given as the total space of an $S^1$-bundle over the total space of a projective bundle over $\C P^1$, where the bundles depend on three parameters $s,t,c\in \Z $ (see Definition \ref{simply connected definition}). The case of a trivial projective bundle (the untwisted case) includes Witten manifolds and was considered in \cite{Wi81,WZ90,KS93}.

The manifolds we are interested in are obtained from certain twisted (i.e. nontrivial) projective bundles and are of dimension $4k+1$. The twisting is necessary to obtain nontrivial relative $\eta $-invariants (see Section \ref{section eta invariants}), which will be used to distinguish components in the moduli space (see Section \ref{section proof}).

Although the construction involves the parameters $s,t$ and $c$, the integral cohomology ring of these manifolds depends up to isomorphism only on $s$ (see Remark \ref{remark cohomology M} and \cite{OU72}). As in the untwisted case, the manifolds can be described as quotients of $S^3\times S^{4k-1}$ by a free action of $S^1$ (see Lemma \ref{product of spheres lemma}).

We now come to the construction of the aforementioned manifolds. Let $\gamma_1$ denote the canonical complex line bundle over $\C P^1$ and let $pr:S^3\to S^3/S^1$ be the Hopf fibration, where $S^3$ is the sphere of unit quaternions. We will always identify $\C P^1$ with $S^3/S^1$ and identify $\gamma _1$ with the complex line bundle $S^3\times _{\rho_1}\C \to S^3/S^1=\C P^1$ associated to the Hopf fibration and the standard one-dimensional representation $\rho _1$ of $S^1$. Similarly, we will identify the $n^{\text {th}}$ tensor power $\gamma _1^n$ with $S^3\times _{\rho _n}\C \to \C P^1$ for $n\in \Z $, where $\rho _n(\lambda)$ for $\lambda \in S^1$ acts on $\C $ via multiplication with $\lambda ^n$.
Let $\epsilon _l$ denote the trivial complex vector bundle of rank $l$ over $\C P^1$. We consider the standard inner product on $\C $ and equip the line bundles above with the induced inner products. For $k>0$ fixed and $c\in \Z $, let $E_c \to \C P^1$ denote the direct sum of $\gamma_1^c$ and $\epsilon _{2k-1}$.

Next we consider the pullback of the bundles above via the projection\linebreak $pr:S^3\to S^3/S^1$. Note that for every $c\in \Z$ there is a canonical trivialization of the complex line bundle $pr^*(\gamma _1^c)$.
 Hence, $pr^*(E_c)\to S^3$ and its associated sphere bundle $S(pr^*(E_c))\to S^3$ have a canonical trivialization. In the following, we will identify $S(pr^*(E_c))$ with $S^3\times S^{4k-1}$ via the corresponding diffeomorphism.

Let $B_c$ be the total space of the projective bundle associated to $E_c\to \C P^1$ and let $q:B_c\to \C P ^1$ denote the projection. Under the identification above $pr^*(B_c)$ corresponds to the quotient of $S^3\times S^{4k-1}$ by $S^1$, where $S^1$ acts trivially on $S^3$ and acts by complex multiplication on $S^{4k-1}\subset \C ^{2k}$. The following lemma follows directly from the description above. The proof is left to the reader.

\begin{lemma}\label{B_c  T^2 lemma}
$B_c$ is diffeomorphic to the quotient of $S^3\times S^{4k-1}$ by a two dimensional torus $T^2$, where $T^2$ acts freely on $S^3\times S^{4k-1}$ by
$$(\lambda,\mu)(x,y):=(x\cdot  \lambda ^{-1}, (\lambda ^c\cdot  y_1\cdot \mu,y_2\cdot \mu, \ldots ,y_{2k}\cdot \mu ))$$
for $(\lambda,\mu) \in T^2$, $ x\in S^3,\; y_1,\ldots ,y_{2k}\in \C \text{ and } y=(y_1,\ldots ,y_{2k})\in S^{4k-1}\subset \C ^{2k}$.\proofend\end{lemma}

Let $u\in H^2(B_c;\Z )$ be the negative of the first Chern class of the canonical complex line bundle over the projective bundle $B_c\to \C P^1$ and let $v$ be the generator of $H^2(\C P^1;\Z )$ defined by $v:=-c_1(\gamma_1)$. 

The tangent bundle along the fibers of $q$, denoted by $T^\triangle$, is a complex vector bundle of rank $2k-1$ over $B_c$ and the tangent bundle $TB_c$ is isomorphic to the complex vector bundle $q^*(T\C P^1)\oplus T^\triangle$ of rank $2k$ (see \cite{BoHi58}). We equip $B_c$ with the induced orientation.

\begin{lemma}\label{B_c topology lemma}
$B_c$ is a simply connected closed oriented $4k$-dimensional manifold. The integral cohomology of $B_c$ as an $H^*(\C P^1;\Z )$-module is given by
$$H^*(B_c;\Z )\cong\Z [u,v]/(v^2,u^{2k}-c\cdot u^{2k-1}\cdot v).$$
In particular, $H^2(B_c;\Z )\cong \Z \langle u,v\rangle $.
Under this identification
$$c(TB_c)=(1+2v)\cdot ((1+u)^{2k}-c\cdot v\cdot (1+u)^{2k-1})$$ and $c_1(TB_c)=(-c+2)\cdot v+2k\cdot u$.
\end{lemma} 

\bigskip
\noindent
{\bf Proof:} Using the homotopy long exact sequence it follows directly that $B_c$ is simply connected. 

By the Leray-Hirsch Theorem, $H^*(B_c;\Z )$ is generated as a $H^*(\C P^1;\Z )$-module by $u$ subject to the relation $u^{2k} +c_1(E_c)\cdot u^{2k-1}+\ldots + c_{2k}(E_c)=0$. Since $c(E_c)=c(\gamma _1 ^{c})=1-c\cdot v$, this gives the statement on the cohomology of $B_c$.

The total Chern class of $T^\triangle$ satisfies $c(T^\triangle)=\sum _{i=0}^{2k} (1+u)^{2k-i}\cdot c_i(E_c)$ (see \cite[p. 514]{BoHi58}). Since $TB_c\cong q^*(T\C P^1)\oplus T^\triangle$, the total Chern class $c(B_c)$ is as stated.
\proofend

\begin{defn}\label{simply connected definition} Let $\overline M_{s,t,c}$ be the total space of the principal $S^1$-bundle over $B_c$ with Euler class equal to $e:=su +tv$. Let $\pi  : \overline M_{s,t,c}\to B_c$ denote the projection.
\end{defn}

\medskip
\noindent
{\bf From now on we will assume that $\mathbf {c}$ is odd, $\mathbf {k\geq 2}$, $\mathbf {s}$ and $\mathbf {t}$ are nonzero coprime integers, and $\mathbf {s}$ is even.}

\begin{lemma}\label{M cohomology lemma}
\begin{enumerate}
\item $\overline M_{s,t,c}$ is simply connected.
\item $H^2(\overline M_{s,t,c};\Z )\cong \Z $, $H^{2i}(\overline M _{s,t,c};\Z )\cong \Z _{s^2}$ is generated by $\pi^*(u)^i$ for\linebreak $4\leq 2i\leq 4k-2$ and $H^{2i+1}(\overline M _{s,t,c};\Z )=0$ for $1\leq 2i+1\leq 4k-3$.
\item $H^{4k-1}(\overline M _{s,t,c};\Z )\cong \Z $, $H^{4k}(\overline M _{s,t,c};\Z )=0$ and $H^{4k+1}(\overline M _{s,t,c};\Z )\cong \Z $.
\item $H^*(\overline M_{s,t,c};\Q)\cong H^*(\C P ^1\times S^{4k-1};\Q )$.\end{enumerate}
\end{lemma}

\bigskip
\noindent
{\bf Proof:} First note that the Euler class $e$ is part of a basis of $H^2(B_c;\Z )\cong \Z ^2$ since $s$ and $t$ are coprime. Using the Gysin sequence for $\overline M_{s,t,c}\to B_c$, one finds that $H^1(\overline M_{s,t,c};\Z )=0$ and $H^2(\overline M_{s,t,c};\Z )\cong \Z$. Hence, $\pi _1(\overline M_{s,t,c})$ vanishes by the Hurewicz theorem and the universal coefficient theorem.

Next note that the cokernel of $\Z \langle v\cdot u^{l-1},u^l\rangle \overset {e\cup}\longrightarrow \Z \langle v\cdot u^{l},u^{l+1}\rangle $ is cyclic of order $s^2$ and generated by $u^{l+1}$ for $1\leq l \leq 2k-2$. The remaining statements now follow from Lemma \ref{B_c topology lemma} and the Gysin sequence.\proofend

\begin{remark}\label{remark cohomology M} Using Poincar\'e duality one finds that the isomorphism type of the ring $H^{*}(\overline M _{s,t,c};\Z )$ depends up to finite ambiguity only on $s$. A closer look shows that $\pi ^*(f)^i$ generates $H^{2i}(\overline M _{s,t,c};\Z )$, $2\leq 2i \leq 4k-2$, where $f$ is chosen so that $(e,f)$ is a basis of $H^2(B_c;\Z )$. It follows that the isomorphism type of the integral cohomology ring of $\overline M _{s,t,c}$ is uniquely determined by $s$.
\end{remark}

Next we consider the smooth $2$-connected cover of $B_c$. Since $H^2(B_c;\Z )\cong \Z ^2$ it can be described as the total space of a principal $T^2$-bundle over $B_c$. We note that the $2$-connected cover is unique up to diffeomorphism

\begin{lemma}\label{product of spheres lemma}
\begin{enumerate}
\item The $2$-connected cover of $B_c$ is diffeomorphic to $S^3\times S^{4k-1}$.
\item $\overline M_{s,t,c}$ is diffeomorphic to a quotient of $S^3\times S^{4k-1}$ by a free action of a subgroup $S^1\subset T^2$.
\end{enumerate}
\end{lemma}

\noindent
{\bf Proof:} The first statement follows directly from Lemma \ref{B_c  T^2 lemma}. For the second statement recall that $s$ and $t$ are coprime. Hence, there is a principal $S^1$-bundle $S^\prime \to B_c$  such that the Euler classes of $\overline M_{s,t,c} \to B_c$ and $S^\prime \to B_c$ generate $H^2(B_c;\Z )$. The two bundles define a principal $T^2$-bundle over $B_c$ with $2$-con\-nected total space. Hence, the latter can be identified with $S^3\times S^{4k-1}$ and $\overline M_{s,t,c}$ is diffeomorphic to a quotient of $S^3\times S^{4k-1}$ by a free action of $S^1$.
\proofend

\begin{remarks}\begin{enumerate}
\item From the homotopy long exact sequence for the fibration
$$S^3\times S^{4k-1}\to \overline M_{s,t,c}$$ one gets $\pi _i(\overline M_{s,t,c})\cong \pi _i(S^3\times S^{4k-1})$ for $i\geq 3$.
\item The principal $T^2$-action on the $2$-connected cover is not equivalent to the standard $T^2$-action on $S^3\times S^{4k-1}$ given by componentwise multiplication since $B_c$ is not diffeomorphic to $\C P^1\times \C P^{2k-1}$.
\item In \cite{OU72}, the cohomology rings of certain $S^1$-quotients of a product of spheres, including  $\overline M_{s,t,c}$, have been computed. The results there can also be used to prove Lemma \ref{M cohomology lemma} and Remark \ref{remark cohomology M}.
\end{enumerate}
\end{remarks}

\section{Diffeomorphism finiteness of $\mathbf {\Z  _2}$-quotients}\label{section diffeo finite}
In this section we show that certain infinite families of $\Z _2$-quotients of the manifolds $\overline M_{s,t,c}$ fall into finitely many oriented diffeomorphism types. {\bf Throughout this section, $\mathbf{s}$ will be a fixed nonzero even integer.}

As before, let $\overline M_{s,t,c}$ be the total space of the principal $S^1$-bundle over $B_c$ with Euler class equal to $su +tv$, where $B_c$ is the total space of the projective bundle associated to $\gamma_1^c\oplus \epsilon _{2k-1}$ and assume that $k\geq 2$, $c$ is odd and $s$ and $t$ are coprime. Let $L\to B_c$ denote the complex line bundle which is associated to the principal $S^1$-bundle.

Consider the total space $M_{s,t,c}$ of the principal $S^1$-bundle over $B_c$ associated to $L\otimes L\to B_c$. Note that the Euler class of $M_{s,t,c}\to B_c$ is equal to $2(su +tv)$.

By construction, $S^1$ acts (fiberwise) on $L$,  $\overline M_{s,t,c}$, $L\otimes L$ and $M_{s,t,c}$. Let $\tau $ denote multiplication by $-1\in S^1$ on the fibers of $L\to B_c$ and on the fibers of $\overline M_{s,t,c}\to B_c$.  Note that the map $L\to L\otimes L$, $v\mapsto v\otimes v$, is equivariant \wrt \ the $\Z _2$-action via $\tau $ on $L$ and the trivial $\Z _2$-action on $L\otimes L$. By passing to the associated principal $S^1$-bundles it follows that $M_{s,t,c}$ can be identified with $\overline M_{s,t,c}/\tau $ and that the quotient map $p:\overline M_{s,t,c}\to M_{s,t,c}$ is a universal covering map.

Since the action of $\tau $ on $\overline M_{s,t,c}$ extends to an action of $S^1$, the fundamental group $\pi _1(M_{s,t,c})=\Z _2$ acts trivially on $\pi _*(\overline M_{s,t,c})$. Hence, $M_{s,t,c}$ is a simple space. In addition, $H^*(M_{s,t,c};\Q )\cong H^*(\overline M_{s,t,c};\Q)^{\pi _1(M_{s,t,c})}\cong H^*(\overline M_{s,t,c};\Q)$, which is isomorphic to $H^*(\C P^1\times S^{4k-1};\Q )$ by Lemma \ref{M cohomology lemma}.

We equip $\overline M_{s,t,c}$ and $M_{s,t,c}$ with the orientation induced from the orientation of $B_c$ (see Section \ref{section family}) and the complex structure of the complex line bundles.

Our aim is to show diffeomorphism finiteness for the family of $(4k+1)$-dimensional oriented manifolds ${\cal F}_s:=\{ M_{s,t,c}\, \mid \, c \text{ and }t \text{ odd and }t \text{ coprime to } s \}$.
\begin{theorem}\label{theorem finiteness} The family ${\cal F}_s$ contains only finitely many oriented diffeomorphism types.
\end{theorem}

\bigskip
\noindent
{\bf Proof:} We first show homotopy finiteness and then diffeomorphism finiteness.

\noindent
{\em Homotopy finiteness claim:} We claim that the family ${\cal F}_s$ belongs to only finitely many simple homotopy types. Note that this is equivalent to showing finiteness of homotopy types since $\pi _1(M_{s,t,c})=\Z _2$ and the Whitehead group of $\Z _2$ is trivial.

Since the members of ${\cal F}_s$ are simple spaces they can be described by Postnikov towers which are classified by their respective $k$-invariants (see for example \cite[Theorem 4.11]{W78}). To show the claim it suffices to prove that there are up to homotopy only finitely many Postnikov towers for the manifolds in this family. 
Let
$$\xymatrix{     &  M_{s,t,c} \ar@{.>}[dl]  \ar[d] \ar[dr] \ar@{.>}[drr] \ar[drrr] \ar[drrrr]&   &  &  & \\
\cdots \ar[r]  &   X_l \ar[r] & X_{l-1} \ar[r] & \cdots \ar[r] & X_1\ar[r] & X_0 }$$
be the Postnikov tower of $M_{s,t,c}$. Recall that each $X_l\to X_{l-1}$ is a principal fibration (with fiber an Eilenberg-MacLane space) which can be described as the pullback of the path fibration
$$K(\pi _{l}(M_{s,t,c}), l)\hookrightarrow \Gamma K(\pi _{l}(M_{s,t,c}), l+1)\to K(\pi _{l}(M_{s,t,c}), l+1)$$ via
a map $\kappa_{l+1}:X_{l-1}\to K(\pi _{l}(M_{s,t,c}), l+1)$. Up to homotopy, the fibration is classified by the homotopy class of $\kappa_{l+1}$, which corresponds to a class $k_{l+1}\in H^{l+1}(X_{l-1};\pi _{l}(M_{s,t,c}))$. As noted before, the Postnikov tower is determined by its $k$-invariants $k_{l}$ for $l\geq 1$. For showing homotopy finiteness, it therefore suffices to show finiteness of the possible $k$-invariants.

By Lemma \ref{product of spheres lemma}, $\overline M_{s,t,c}$ is a quotient of $S^3\times S^{4k-1}$ by a free action of $S^1$. Since $\overline M_{s,t,c}$ is simply connected and $M_{s,t,c}$ is the quotient of a free $\Z _2$-action on $\overline M_{s,t,c}$, we have $\pi _1( M_{s,t,c})=\Z _2$, $\pi_2(M_{s,t,c})\cong \pi_2(\overline M_{s,t,c})\cong \Z $ and $\pi _l(M_{s,t,c})\cong \pi _l(S^3\times S^{4k-1})$ for $l\geq 3$.

It follows that the stages $X_{\leq 2}$ of the Postnikov tower of $M_{s,t,c}$ do not depend, up to homotopy, on the choice of the parameters. In fact, one has $X_0=\{pt\}$,
$$\kappa_2:X_0\to K(\Z _2,2), \, k_2=0, \, X_1\simeq K(\Z _2,1)\simeq \R P^\infty \text{ and}$$
$$\kappa_3:X_1\to K(\Z ,3), \, k_3\in H^3(\R P^\infty ;\Z )=0, \, X_2\simeq X_1\times K(\Z ,2)\simeq \R P^\infty \times \C P^\infty .$$

In the following we will consider the stages $X_l$ and partial Postnikov towers up to homotopy without explicit mention. 

The next stage $X_3$ in the Postnikov tower is determined by the invariant $$k_4\in H^4(X_2; \pi _3(S^3\times S^{4k-1}))\cong H^4(\R P^\infty \times \C P^\infty;\Z )\cong \Z \times \Z _2 \times \Z _2$$ (recall that $k\geq 2$). Using the Gysin sequence, one finds that $\vert H^4(M_{s,t,c};\Z)\vert =4s^2$. Since $X_3$ is obtained from $M_{s,t,c}$ by attaching cells of dimension $\geq 5$, the homomorphism $H^4(X_3;\Z )\to H^4(M_{s,t,c};\Z )$ is injective. Hence, the cohomology group $H^4(X_3;\Z )$ is finite and determined up to finite ambiguity by $s$. The invariant $k_4$  can be identified with the transgression of the fundamental class of the fiber in the Leray-Serre spectral sequence for the fibration $X_3\to X_2$. It follows that $k_4$ is determined up to finite ambiguity by $s$. Hence, $X_3$ is determined up to finite ambiguity by $s$ as well. For later reference we note that $H^{\geq 3}(X_3;\Q )=0$ since $s\neq 0$ (again by applying the Leray-Serre spectral sequence).

Since $\pi _l(S^3\times S^{4k-1})\otimes \Q =0$ for $3<l<4k-1$, it follows by induction that, for $l< 4k-1$, the invariants $k_{l+1}$ for $M_{s,t,c}$ and its stages $X_l$ are determined, up to finite ambiguity, by $s$. Hence, the same holds for the partial Postnikov tower $(X_{4k-2}\to X_{4k-3}\to \ldots \to X_1\to X_0)$ of $M_{s,t,c}$. Again by induction, or by using the minimal model, one finds that $H^{\geq 3}(X_l;\Q )=0$ for $2<l<4k-1$.

By the above the invariant $k_{4k}\in H^{4k}(X_{4k-2};\pi _{4k-1}(M_{s,t,c}) )$ is also determined up to finite ambiguity by $s$. Since $\pi _{>4k-1}(M_{s,t,c})\otimes \Q =0$ we can argue as before to see that, for every $l\geq 4k-1$, the partial Postnikov tower $(X_{l}\to X_{l-1}\to \ldots \to X_1\to X_0)$ of $M_{s,t,c}$ is also determined up to finite ambiguity by $s$.  Moreover, the construction of the infinitely many stages $X_l$ for $l{>4k+1}$ of the Postnikov tower is formal, i.e. it depends only on $X_{4k+1}$ (see \cite[p. 72]{GM13}). Hence, the entire tower is determined up to finite ambiguity by $s$, and the claim follows.

\medskip
\noindent
{\em Diffeomorphism finiteness claim:} We claim that after restricting to a (simple) homotopy type there are only finitely many oriented diffeomorphism types among the manifolds $M_{s,t,c}$.
Let us fix a homotopy type represented by $M\in {\cal F}_s$ and consider the subfamily ${\cal F}_s^\prime :=\{M_{s,t,c}\in {\cal F}_s\, \mid \, M_{s,t,c}\simeq M\}$ of manifolds homotopy equivalent to $M_{s,t,c}$. Recall that each $M_{s,t,c}$ comes with an orientation. To show that the family ${\cal F}_s^\prime $ contains only finitely many oriented diffeomorphism types, we apply the surgery exact sequence \cite{W70}
$$\cdots \to L_{4k+2}(\Z _2)\to {\cal S}(M)\to [M,G/O]\to \cdots .$$
Note that $H^*(M;\Q )\cong H^*(\C P^1\times S^{4k-1};\Q )$ and the homotopy groups $\pi _i(G/O)$ of the $H$-space $G/O$ are finite for $i\not \equiv 0 \bmod 4$. Hence, $[M,G/O]$ is finite. Since $L_{4k+2}(\Z _2) = \Z _2$, the smooth structure set ${\cal S}(M)$ is also finite, and the claim follows.

Combining the two claims above, we conclude that for fixed $s$ there are, up to orientation preserving diffeomorphism, only finitely many $(4k+1)$-dimensional manifolds in the family ${\cal F}_s$.
\proofend

\section{Nonnegative sectional and positive Ricci curvature}\label{section nonnegative curvature}
In this section we consider submersion metrics of nonnegative sectional and positive Ricci curvature on $M_{s,t,c}$ and $\overline M_{s,t,c}$ and extend the latter to the associated disk bundle. 

Let $(S^l,h_{S^l})$ denote the round sphere of radius $1$ and let $h_{S^3}\times h_{S^{4k-1}}$ denote the product metric on $S^3\times S^{4k-1}$. Recall from Lemma \ref{B_c  T^2 lemma} that $T^2$ acts freely and isometrically on $(S^3\times S^{4k-1},h_{S^3}\times h_{S^{4k-1}})$ with quotient diffeomorphic to $B_c$. By Lemma \ref{product of spheres lemma}, $\overline M_{s,t,c}$ is diffeomorphic to a quotient of $S^3\times S^{4k-1}$ by an $S^1$-subaction of $T^2$. Let $\overline g_{s,t,c}$ denote the submersion metric on $\overline M_{s,t,c}$, i.e. $(S^3\times S^{4k-1},h_{S^3}\times h_{S^{4k-1}})\to (\overline M_{s,t,c},\overline g_{s,t,c})$ is a Riemannian submersion. We note that $M_{s,t,c}$ can be identified with the quotient of $S^3\times S^{4k-1}$ by a subgroup of $T^2$ which is isomorphic to $S^1\times \Z _2$. Let $g_{s,t,c}$ denote the submersion metric on $M_{s,t,c}$. By construction $p :(\overline M_{s,t,c},\overline g_{s,t,c})\to (M_{s,t,c},g_{s,t,c})$ is a Riemannian universal covering.

\begin{lemma}\label{curvature lemma} $(M_{s,t,c},g_{s,t,c})$ and
$(\overline M_{s,t,c},\overline g_{s,t,c})$ both have nonnegative sectional and positive Ricci curvature.
\end{lemma}

\bigskip
\noindent
{\bf Proof:} Note that the sectional curvature of $(S^3\times S^{4k-1},h_{S^3}\times h_{S^{4k-1}})$ is always nonnegative and vanishes only on mixed planes. It is easy to see that there is, for any horizontal vector of the Riemannian submersion $S^3\times S^{4k-1}\to M_{s,t,c}$ (resp. $S^3\times S^{4k-1}\to \overline M_{s,t,c}$), a horizontal plane of positive sectional curvature which contains this vector. Hence, the statements follow from the Gray-O'Neill formula  \cite{Gr67,ON66}.\proofend

Recall that $\overline M_{s,t,c}$ (resp. $M_{s,t,c}$) is a quotient of $S^3\times S^{4k-1}$ by a subgroup $H\subset T^2$ which is isomorphic to $S^1$ (resp. $S^1\times \Z _2$). We remark that the normalizer $N$ of $H$ in the isometry group of $(S^3\times S^{4k-1},h_{S^3}\times h_{S^{4k-1}})$ acts with cohomogeneity $1$ on $\overline M_{s,t,c}$ (resp. $M_{s,t,c}$) and the metric $\overline g_{s,t,c}$ (resp. $g_{s,t,c}$) is $N$-invariant.

For the computation of $\eta$-invariants in the next sections, we will also need to put a suitable metric on the disk bundle associated to the principal $S^1$-bundle $\overline M_{s,t,c}\to B_c$. Let $W_{s,t,c}:=\overline M_{s,t,c}\times _{S^1}D^2$, where $D^2\subset \R ^2$ is the disk of radius one. We equip $D^2$ with a metric $g_{D^2}$ (a torpedo metric) such that $g_{D^2}$ is $S^1$-invariant, is of product type on the annulus $\{ x\in D^2\, \mid \, \vert x\vert \geq 1-\epsilon \}$ for a fixed small positive $\epsilon  $, and such that $g_{D^2}$ is of positive curvature outside of the annulus. Next we consider the product metric $\overline g_{s,t,c}\times g_{D^2}$ on $\overline M_{s,t,c} \times D^2$ and denote by $h_{s,t,c}$ the submersion metric on $W_{s,t,c}$ \wrt \ the quotient map $\overline M_{s,t,c} \times D^2 \to W_{s,t,c}$. The next lemma follows directly from the construction and the Gray-O'Neill formula \cite{Gr67,ON66}.
\begin{lemma} The metric $h_{s,t,c}$ extends $\overline g_{s,t,c}$ to an $S^1$-invariant metric on $W_{s,t,c}$ of nonnegative sectional and positive scalar curvature which is of product type near the boundary.\proofend
\end{lemma}

\section{$\mathbf{Spin^\ccc}$-structures and Dirac operators}\label{section spinc structures}
In this section we introduce suitable $Spin^\ccc$-structures and corresponding Dirac operators on $(M_{s,t,c},g_{s,t,c})$, on its universal cover and on the associated disk bundle. These will be used to compute $\eta$-invariants in the next section. For background information on and references for $Spin^\ccc$-manifolds and Dirac operators we refer to \cite{ABS64,APSI75,LM89,DG21}.

We begin by defining the relevant $Spin^\ccc$-structures. Recall that $\pi $ denotes the projection $\overline M_{s,t,c}\to B_c$. In the following we will also denote the projections $W_{s,t,c}\to B_c$ and $M_{s,t,c}\to B_c$ by $\pi$. Also we will suppress the parameters $s,t$ and $c$ in the notation for $Spin^\ccc$-structures and Dirac operators.

Recall that $\tau $ acts freely by multiplication with $-1\in S^1$ on the fibers of $\overline M_{s,t,c}\to B_c$ and that the quotient can be identified with $M_{s,t,c}$. Let $\tau $ also denote the action by $-1$ on the fibers of the disk bundle $W_{s,t,c}\to B_c$.

The action of $\Z _2=\{\mathrm{id}, \tau \}$ on $\overline M_{s,t,c}$ and $W_{s,t,c}$ lifts via differentials to the respective oriented orthonormal frame bundles.

\begin{lemma}\label{lemma unique spin structure}
$(W_{s,t,c},h_{s,t,c})$ has a unique $Spin$-structure.
\end{lemma}

\bigskip
\noindent
{\bf Proof:} Recall that $c$ and $t$ are odd and $s$ is even. Since $TW_{s,t,c}\cong \pi ^*(TB_c\oplus L)$, $c_1(TB_c)=(-c+2)\cdot v+2k\cdot u$ and $c_1(L)=su +tv$ (see Section \ref{section family}), the manifold $W_{s,t,c}$ is spin. Moreover, the $Spin$-structure on $(W_{s,t,c},h_{s,t,c})$ is unique since $H^1(W_{s,t,c};\Z _2)=0$.\proofend

We note that the induced structure on the boundary is the unique $Spin$-structure on $(\overline M_{s,t,c},\overline g_{s,t,c})$ since $\pi_1(\overline M_{s,t,c})=0$. Note, however, that $M_{s,t,c}$ is not spin but does admit a $Spin^\ccc$-structure, as will be explained below. It also carries a twisted $Spin$-structure in the sense of \cite{BG96}.

Let $P_{SO}(W)\to W_{s,t,c}$ be the principal bundle of oriented orthonormal frames and let $P_{Spin}(W)\to P_{SO}(W)$ be the covering map defining the $Spin$-structure. Its restriction to a fiber of $P_{Spin}(W)\to W_{s,t,c}$ can be identified (noncanonically) with the non-trivial covering $\rho: Spin(4k+2)\to SO(4k+2)$.

The fixed point manifold of the $\tau$-action on $W_{s,t,c}$ is the zero section $B_c$, which is of codimension $2$. Hence, the involution $\tau $ is of odd type and the $\Z _2$-action on $P_{SO}(W)$ does not lift to the $Spin$-structure (see \cite[p. 487]{AB68}). However, as we will see below, the $\Z _2$-action does lift to a suitable $Spin^\ccc$-structure.

Let $P_{U(1)}(W)\to W_{s,t,c}$ be the trivial principal $U(1)$-bundle and consider the two-fold covering map $P_{U(1)}(W)\to P_{U(1)}(W)$ for which the restriction to a fiber is given by the non-trivial two-fold covering $ (\;\;  )^2: U(1)\to U(1)$, $\lambda \mapsto \lambda ^2$.

Let $\Z _2$ act by multiplication with $\pm 1 $ on $U(1)$. The $\Z _2$-actions on $W_{s,t,c}$ and $U(1)$ define a $\Z _2$-action on $P_{U(1)}(W)$. Note that this $\Z _2$-action does not lift in the two-fold covering $P_{U(1)}(W)\to P_{U(1)}(W)$.

Let $P_{Spin^\ccc }(W)\to W_{s,t,c}$ denote the $Spin^\ccc$-structure associated to the $Spin$-structure on $W_{s,t,c}$.

\begin{lemma}
The $\Z _2$-actions on $P_{SO}(W)$ and $P_{U(1)}(W)$ lift to a $\Z _2$-action on $P_{Spin^\ccc }(W)$.
\end{lemma}

\bigskip
\noindent
{\bf Proof:} By definition, $P_{Spin^\ccc }(W)$ is the extension of $P_{Spin}(W)$ \wrt \ the inclusion
$$Spin(4k+2)\hookrightarrow (Spin(4k+2)\times U(1))/\{\pm (1,1) \}=Spin ^\ccc(4k+2).$$
Moreover, there is a $Spin^\ccc (4k+2)$-equivariant bundle map
$$P_{Spin^\ccc }(W) \longrightarrow P_{SO(n)}(W) \times P_{U(1)}(W)$$ \wrt \ 
the homomorphism $Spin^\ccc (4k+2)\xrightarrow {\rho \times (\; \; )^2} SO(4k+2)\times U(1)$ (here $P_{SO(n)}(W) \times P_{U(1)}(W)$ denotes the fiberwise product of $P_{SO(n)}(W)$ and $P_{U(1)}(W)$).

Recall that the $\Z _2$-actions on $P_{SO}(W)$ and $P_{U(1)}(W)$ do not lift as $\Z _2$-actions in the coverings $P_{Spin}(W)\to P_{SO}(W)$ and $P_{U(1)}(W)\to P_{U(1)}(W)$. In both cases the induced action on the total spaces is by an effective action of $\Z _4$. Note, however, that the diagonal action of $\Z _4$ on $P_{Spin^\ccc }(W)$ has $\Z _2\subset \Z _4$ as ineffective kernel. Hence, the $\Z _2$-action on $P_{SO}(W)\times P_{U(1)}(W)$ lifts as a $\Z _2$-action to the $Spin^\ccc$-structure  $P_{Spin^\ccc }(W)\to W_{s,t,c}$.\proofend

Recall that the $\Z _2$-actions on $(W_{s,t,c},h_{s,t,c})$ and on the trivial principal $U(1)$-bundle $P_{U(1)}(W)\to W_{s,t,c}$ are of product form near the boundary of $W_{s,t,c}$. We fix a flat unitary $\Z _2$-equivariant connection $\nabla ^\ccc (W)$ on $P_{U(1)}(W)\to W_{s,t,c}$ which is constant in the normal direction near the boundary of $W_{s,t,c}$.

Next we describe the relevant Dirac operators on $W_{s,t,c}$ and its boundary. Let $S(W_{s,t,c})$ denote the spinor bundle for the $\spinc$-structure on $W_{s,t,c}$ defined before.
The Levi-Civita connection of $(W_{s,t,c},h_{s,t,c})$ together with the connection $\nabla ^\ccc (W)$ determine a connection $\nabla(W)$ on $S(W_{s,t,c})$. Let $D_W$ be the associated $\spinc$-Dirac operator, i.e. $D_W$ is the composition
\begin{small}$$
\Gamma (S(W_{s,t,c}))\to \Gamma (S(W_{s,t,c})\otimes T^*W_{s,t,c})\to  \Gamma (S(W_{s,t,c})\otimes TW_{s,t,c})\to \Gamma (S(W_{s,t,c})),$$
\end{small}where the first map is the connection $\nabla(W)$, the second map uses the isomorphism given by the metric $h_{s,t,c}$ and the last map is induced from Clifford multiplication  (see \cite[D.9]{LM89}).

Since $W_{s,t,c}$ is of even dimension the spinor bundle $S(W_{s,t,c})$ splits as a direct sum $S^+(W_{s,t,c})\oplus S^-(W_{s,t,c})$ and the operator $D_W$ restricts to an operator
$${\mathfrak D_W}^+:\Gamma (S^+(W_{s,t,c}))\to \Gamma (S^-(W_{s,t,c})).$$

The $\spinc$-structure on $W_{s,t,c}$ induces a $\spinc$-structure on the boundary. Let $\overline P\to \overline M_{s,t,c}$ denote the corresponding principal $\spinc$-bundle. The restriction of $S^+(W_{s,t,c})$  and ${\mathfrak D_W}^+$ to the boundary can be identified with the spinor bundle $S(\overline M_{s,t,c})$ and the $\spinc$-Dirac operator
$$D_{\overline M}:\Gamma (S(\overline M_{s,t,c}))\to \Gamma (S(\overline M_{s,t,c}))$$ on $(\overline M_{s,t,c},\overline g_{s,t,c})$, which is defined \wrt \  $\overline P\to \overline M_{s,t,c}$ and the restriction $\overline \nabla ^\ccc$ of the connection  $\nabla ^\ccc (W)$ to the principal $U(1)$-bundle (see \cite{APSI75}) $$\overline P _{U(1)}:=P_{U(1)}(W)\vert_{\overline M_{s,t,c}}\to \overline M_{s,t,c}.$$
Consider the orthogonal projection of $\Gamma (S^+(W_{s,t,c})\vert _{\overline M_{s,t,c}})=\Gamma (S(\overline M_{s,t,c}))$ onto the space spanned by the eigenfunctions of $D_{\overline M}$ for nonnegative eigenvalues. Following Atiyah, Patodi and Singer we impose the \emph{APS-boundary condition}, i.e. we  restrict to sections $\phi\in \Gamma(S^+(W_{s,t,c}))$ for which $\phi\vert_{\overline M_{s,t,c}}$ is in the kernel of the projection. After imposing this condition, the operator ${\mathfrak D_{W}}^+$ has finite dimensional kernel and will be denoted by $D_W^+$. Similarly, the formal adjoint of ${\mathfrak D_W}^+$ (defined via bundle metrics) subject to the adjoint APS-boundary condition has finite dimensional kernel and will be denoted by $(D_W^+)^*$. The index of $D_W^+$ is defined as $\mathrm{ind} \, D_W^+:= \dim\ker \, D_W^+ - \dim\ker \,  (D_W^+)^* \in \Z $ (see \cite{APSI75} for details). Note that by construction the operators $D_W^+$, $(D_W^+)^*$ and $D_{\overline M}$ are $\Z _2$-equivariant. For later reference we point out the following crucial lemma:
\begin{lemma}\label{psc index lemma} The operators $D_W^+$,  $(D_W^+)^*$ and $D_{\overline M}$ are injective. In particular, we have $\mathrm{ind} \, D_W^+=0$.
\end{lemma}

\bigskip
\noindent
{\bf Proof:} Since $h_{s,t,c}$ and $\overline g_{s,t,c}$ are of positive scalar curvature and all relevant connections are flat, the statements follow from the argument of Schr\"odinger and Lichnerowicz \cite{S32,Li63,LM89}.\proofend

Note that the objects on $W_{s,t,c}$ considered above, when restricted to the boundary $\overline M_{s,t,c}$, induce corresponding objects on $M_{s,t,c}$ by passing to the quotient \wrt \ the $\Z _2$-action. For example, the quotient of $(\overline M_{s,t,c},\overline g_{s,t,c})$ by the free  isometric $\Z _2$-action can be identified with $(M_{s,t,c},g_{s,t,c})$ and the same is true for the respective principal bundles of oriented orthonormal frames and the Levi-Civita connections.

Similarly, the $\Z _2$-quotient of the principal $U(1)$-bundle $\overline P _{U(1)}\to \overline M_{s,t,c}$ with its flat connection $\overline \nabla ^\ccc$ and the quotient of the $Spin^\ccc$-structure $\overline P\to \overline M_{s,t,c}$ can be identified with a principal $U(1)$-bundle $P _{U(1)}\to M_{s,t,c}$ with flat connection $\nabla ^\ccc$ and a $Spin^\ccc$-structure $P\to M_{s,t,c}$ on $M_{s,t,c}$, respectively.

Since the generator of $\Z _2$ acts  by $(\tau ,-1)$ on $\overline P_{U(1)}=\overline M_{s,t,c}\times U(1)$ the bundle $P _{U(1)}\to M_{s,t,c}$ can be identified with $\overline M_{s,t,c}\times _{\Z _2}U(1)\to M_{s,t,c}$. This bundle is nontrivial. In fact, its first Chern class is of order $2$ and generates the kernel of $p^*: H^2(M_{s,t,c};\Z )\to H^2(\overline M_{s,t,c};\Z )$.

Let $S(M_{s,t,c})$ denote the spinor bundle associated to the $Spin^\ccc$-structure on $M_{s,t,c}$ and let
$$D_{M}:\Gamma (S(M_{s,t,c}))\to \Gamma (S(M_{s,t,c}))$$
denote the associated $\spinc$-Dirac operator. It follows from the construction that $D_{M}$ lifts to the 
$\Z _2$-equivariant $\spinc$-Dirac operator $D_{\overline M}$ \wrt \ the covering map $p:\overline M_{s,t,c}\to M_{s,t,c}$.

\section{Computation of $\mathbf{\eta}$-invariants}\label{section eta invariants}
In this section we will compute relative $\eta$-invariants for the $Spin^\ccc$-Dirac operator $D_{M}$ on $M_{s,t,c}$ twisted with the nontrivial complex $1$-dimensional representation of $\pi _1(M_{s,t,c})$.
These computations will be used in the next section to prove the main theorem. Alternatively, one could use a twisted $Spin$-structure of $M_{s,t,c}$ and compute $\eta$-invariants of associated Dirac operators along the lines of \cite[Section 2]{BG96}.

The idea to use relative $\eta$-invariants to distinguish components of moduli spaces goes back to Atiyah, Patodi and Singer who explained this for positive scalar curvature metrics on spin manifolds in \cite{APSII75}. They also pointed out the possibility to extend this idea to certain $Spin^\ccc$-manifolds. For background information on $\eta$-invariants of $Spin^\ccc$-manifolds, we also refer to \cite{DG21}. 

Recall that $\pi _1(M_{s,t,c})=\Z _2$ and $p:\overline M_{s,t,c}\to M_{s,t,c}$ is a universal covering. Let $\alpha :\pi _1(M_{s,t,c})\to U(1)$ denote the nontrivial homomorphism and let $\alpha $ also denote the associated complex line bundle  ${\overline M_{s,t,c}}\times _\alpha \C \to M_{s,t,c}$. We fix a flat unitary connection on $\alpha $. Let $D_{M,\alpha }$ denote the $Spin^\ccc$-Dirac operator $D_{M}$ twisted with $\alpha$.

Next consider the $\eta$-invariants $\eta (M_{s,t,c})$ and $\eta _\alpha  (M_{s,t,c})$ of $D_{M}$ and $D_{M,\alpha }$, respectively. Recall that $\eta (M_{s,t,c})$ (resp. $\eta _\alpha  (M_{s,t,c})$) is given by the value at $z=0$ of the meromorphic extension of the series $\sum_\lambda \frac{\sign(\lambda)}{\vert\lambda \vert^{z}}$ for $z\in\C $ with $\mathrm{Re} (z)\gg 0$
to the complex plane, 
where the sum is taken over all nonzero eigenvalues $\lambda$ of $D_M$ (resp. $D_{M,\alpha }$) (see \cite{APSI75} for background information on $\eta$-invariants).

\begin{defn}
The relative $\eta$-invariant $\tilde \eta _\alpha(M_{s,t,c})$ is defined by
$$\tilde \eta _\alpha (M_{s,t,c}):=\eta _\alpha  (M_{s,t,c}) - \eta (M_{s,t,c}).$$
\end{defn}

Recall from Section \ref{section spinc structures} that $D_{M}$ lifts to the 
$\Z _2$-equivariant $\spinc$-Dirac operator $D_{\overline M}$. The $\eta$-invariant of $D_{\overline M}$ refines to a $\Z _2$-equivariant invariant with values denoted by $\eta (M_{s,t,c})_g$ for $g\in \Z _2=\{1,\tau\}$. As pointed out in \cite[Theorem 3.4]{D78}, the $\eta$-invariants for $M_{s,t,c}$ can be computed from equivariant $\eta$-invariants for $\overline M_{s,t,c}$. In our situation this relation is given by

\begin{align*}
\eta _\alpha (M_{s,t,c}) & = \frac { \eta (\overline M_{s,t,c})_1\cdot \chi_\alpha (1) + \eta (\overline M_{s,t,c})_\tau \cdot \chi_\alpha (\tau)} 2=\frac {\eta (\overline M_{s,t,c})-\eta (\overline M_{s,t,c})_\tau} 2\end{align*}
and
\begin{align*}\eta (M_{s,t,c}) & =\eta _e(M_{s,t,c}) = \frac {\eta  (\overline M_{s,t,c})_1\cdot 1 + \eta  (\overline M_{s,t,c})_\tau \cdot 1} 2 =\frac {\eta (\overline M_{s,t,c}) +\eta  (\overline M_{s,t,c})_\tau } 2,
\end{align*}
where $\chi_\alpha$ is the character of $\alpha $ and $e:\pi _1(M_{s,t,c})\to  U(1)$ denotes the trivial representation. This gives, for the relative $\eta$-invariant,
$$
\widetilde{\eta}_\alpha (M_{s,t,c})=\eta _\alpha (M_{s,t,c})-\eta (M_{s,t,c})=-\eta (\overline M_{s,t,c})_\tau .
$$

Next we consider the $\Z _2$-action on the disk bundle $W_{s,t,c}$ over $B_c$ and the equivariant $Spin^\ccc$-Dirac operator $D_W^+$ which was defined in Section \ref{section spinc structures}. Since $\tau$ acts by $-1$ on the fibers of $W_{s,t,c}$ the fixed point manifold can be identified with $B_c$. Let $a(B_c)(\tau)$ be the local datum of the Lefschetz fixed-point formula for the $\Z _2$-equivariant operator $D_W^+$ at $B_c$ evaluated at $\tau \in\Z _2$ as described in \cite{ASIII68} (see also \cite{DG21}).

The index formula for manifolds with boundary \cite{APSI75} refines in the presence of symmetries and gives a relation between equivariant $\eta$-invariants, local data and certain representations attached to the index of $D_W^+$  and the kernel of $D_{\overline M}$ (see \cite{D78} for details). In our situation one obtains

\begin{proposition}\label{proposition eta and local datum}
$\widetilde{\eta}_\alpha (M_{s,t,c})=-2a(B_c)(\tau )$.
\end{proposition}

\bigskip
\noindent
{\bf Proof:} As a warm-up we first consider the nonequivariant APS-index formula for $D_W^+$ which takes the form (see \cite[Theorem 3.10 and Section 4]{APSI75}) 
$$\mathrm{ind} \, D_W^+  = \left ( \int _{W_{s,t,c}} e^{\frac 1 2 c_1}\hat{{\cal A}}(W_{s,t,c},h_{s,t,c}) \right ) - \frac{\dim h(\overline M_{s,t,c},\overline g_{s,t,c}) + \eta (\overline M_{s,t,c})}{2},$$
where $c_1$ denotes the first Chern form of $\nabla ^\ccc (W)$, $\hat{{\cal A}}(W_{s,t,c},h_{s,t,c})$ represents the $\hat {\cal A}$-series evaluated on the Pontryagin forms $p_i(W_{s,t,c},h_{s,t,c})$ and $h(\overline M_{s,t,c},\overline g_{s,t,c})$ is the kernel of $D_{\overline M}$. Since $\nabla ^\ccc (W)$ is flat, $c_1$ vanishes. Since $(W_{s,t,c},h_{s,t,c})$ and $(\overline M_{s,t,c},\overline  g_{s,t,c})$ are of positive scalar curvature, $\mathrm{ind} \, D_W ^+ $ and $h(\overline M_{s,t,c},\overline g_{s,t,c})$ both vanish (see Lemma \ref{psc index lemma}). 

Next we consider the index of the $\Z _2 $-equivariant operator $D_W^+$. The index evaluated at $\tau $ can be expressed by the formula above after making the following replacements  (see \cite[Theorem 1.2]{D78} for details): First, $\dim h(\overline M_{s,t,c},\overline g_{s,t,c})$ is replaced by the character of the $\Z _2$-representation given by the kernel of $D_{\overline M}$ evaluated at $\tau$. We denote this value by $h_{\tau}$. Next, $\eta (\overline M_{s,t,c})$ is replaced by $\eta (\overline M_{s,t,c})_\tau$. Finally, the integral is replaced by the local datum $a(B_c)(\tau )$. Hence, one has
$$\mathrm{ind} \, D_W^+ (\tau ) = a(B_c)(\tau ) - \frac{h_{\tau} + \eta (\overline M_{s,t,c})_\tau}{2}.$$
 Since $(W_{s,t,c},h_{s,t,c})$ and $(\overline M_{s,t,c},\overline  g_{s,t,c})$ are of positive scalar curvature, the representations which are used to define $\mathrm{ind} \, D_W^+ (\tau ) $ and $h_{\tau}$ are all $0$-dimensional and trivial (see Lemma \ref{psc index lemma}). Hence,   $\mathrm{ind} \, D_W^+ (\tau ) $ and $h_{\tau}$ both vanish and
$$\widetilde{\eta}_\alpha (M_{s,t,c})=-\eta (\overline M_{s,t,c})_\tau = -2a(B_c)(\tau).$$\proofend

We proceed to describe the local datum $a(B_c)(\tau )$ (see \cite[Section 3]{ASIII68} for the general discussion). Let $\{\pm x_1,\ldots , \pm x_{2k}\}$ denote the formal roots of $TB_c$ and let $y$ be the Euler class of the oriented normal bundle $\nu _{B_c}$ of 
$B_c\subset W _{s,t,c}$. Let $c_1$ denote now the first Chern class of the $Spin^\ccc$-Dirac operator $D_W^+$. Then the local datum evaluated at $\tau $ is given by
$$a(B_c)(\tau)=\epsilon \cdot \int _{B_c} e^{\frac{1}{2}c_1}\cdot \hat{{\cal A}}(B_c)\cdot \frac 1 {i\cdot e^{y/2} + i\cdot e^{-y/2}},$$
where $\hat{{\cal A}}(B_c)=\prod _{j=1}^{2k} \frac {x_j}{e^{x_j/2}-e^{-x_j/2}}$ and $\epsilon \in \{\pm i \}$ depends on the lift of the $\Z _2$-action to the $Spin^\ccc$-structure. We will not discuss this ambiguity further since it will not affect the results on moduli spaces stated in the main theorem. The class  $c_1$ vanishes since the bundle $P_{U(1)}(W)$ is trivial. Note that $y=su +tv$ since $\nu _{B_c}$ is isomorphic to the complex line bundle associated to the principal $S^1$-bundle $\pi  : \overline M_{s,t,c}\to B_c$  (see Definition \ref{simply connected definition}).
Hence,
$$a(B_c)(\tau)=\pm \int _{B_c} \hat{{\cal A}}(B_c)\cdot \frac 1 {e^{(su +tv)/2} + e^{-(su +tv)/2}}.$$

Next recall from Lemma \ref{B_c topology lemma} that $TB_c$ has a complex structure and the total Chern class of $TB_c$ is given by
$$c(TB_c)=(1+2v)\cdot ((1+u)^{2k}-c\cdot v\cdot (1+u)^{2k-1})=(1+2v)\cdot (1+u)^{2k-1}\cdot (1+u-c\cdot v).$$ Hence, one obtains the local term $a(B_c)(\tau)$ up to sign by integrating
\begin{equation}\label{local datum eq} \frac {2v}{e^v-e^{-v}}\cdot \left (\frac u {e^{\frac u 2}-e^{-\frac u 2}}\right )^{2k-1}\cdot \frac {u-cv}{e^{\frac {u-cv} 2}-e^{- \frac {u-cv} 2}}\cdot \frac 1 {e^{\frac {su +tv} 2} + e^{- \frac {su +tv} 2}}\end{equation}
over $B_c$. Note that the integral is given by evaluating the cohomological expression on the fundamental class of $B_c$, which, by Lemma \ref{B_c topology lemma}, amounts to computing the coefficient of $u^{2k-1}\cdot v$ in (\ref{local datum eq}). In the following, $k$ and $c$ will be fixed.

\begin{proposition}\label{proposition local datum} For almost all $s\neq 0$ with $s$ even,
$a(B_c)(\tau)$ is a non-zero polynomial in $t$ of degree $1$.
\end{proposition}
The proposition as stated is sufficient for our purposes. It is likely that the statement is true for all $s\neq 0$. We leave it to the interested reader to prove the more general statement.

\bigskip
\noindent
{\bf Proof:} For a fixed odd integer $c$, let $A\in \Q [s,t]$ denote the polynomial obtained by integrating the expression in (\ref{local datum eq}) over $B_c$. To prove the proposition we first note that the factor in (\ref{local datum eq}) involving $t$
 is equal to $$\frac 1 {e^{\frac {su} 2} + e^{- \frac {su} 2}}\cdot \left (1- \frac {tv} 2\cdot \frac {e^{\frac {su} 2} - e^{- \frac {su} 2}}{e^{\frac {su} 2} + e^{- \frac {su} 2}}\right )$$
since $v^2=0$ by Lemma \ref{B_c topology lemma}.

Hence, $A$ is a polynomial in $t$ of degree $\leq 1$, say $A=A_0-A_1\cdot t$ with $A_i\in \Q [s]$. Moreover, by looking at the other factors of (\ref{local datum eq}), we see that $A_1$ is given by integrating
$$\frac {2v}{e^v-e^{-v}}\cdot \left (\frac u {e^{\frac u 2}-e^{-\frac u 2}}\right )^{2k-1}\cdot \frac {u-cv}{e^{\frac {u-cv} 2}-e^{- \frac {u-cv} 2}}\cdot\frac 1 {e^{\frac {su} 2} + e^{- \frac {su} 2}}\cdot \left (\frac {v} 2\cdot \frac {e^{\frac {su} 2} - e^{- \frac {su} 2}}{e^{\frac {su} 2} + e^{- \frac {su} 2}}\right )
$$
over $B_c$. Since $v^2=0$ we get
$$A_1=\int _{B_c} \left (\frac u {e^{\frac u 2}-e^{-\frac u 2}}\right )^{2k-1}\cdot \frac {u}{e^{\frac {u} 2}-e^{- \frac {u} 2}}\cdot\frac 1 {e^{\frac {su} 2} + e^{- \frac {su} 2}}\cdot \frac {e^{\frac {su} 2} - e^{- \frac {su} 2}}{e^{\frac {su} 2} + e^{- \frac {su} 2}}\cdot \frac {v} 2.
$$
Using Lemma \ref{B_c topology lemma} again, it follows that $A_1$ is equal to the coefficient of $u^{2k-1}$ in the formal power series
$$\left (\frac u {e^{\frac u 2}-e^{-\frac u 2}}\right )^{2k}\cdot\frac {e^{\frac {su} 2} - e^{- \frac {su} 2}} {2(e^{\frac {su} 2} + e^{- \frac {su} 2})^2}\in \Q [s][[u]].$$
 Note that $A_1$  is an odd polynomial in $s$ of degree $\leq 2k-1$, which can be written as a residue,
$$A_1=\mathrm {Res}_{u=0}\left (\left (\frac 1 {e^{\frac u 2}-e^{-\frac u 2}}\right )^{2k}\cdot\frac {e^{\frac {su} 2} - e^{- \frac {su} 2}} {2(e^{\frac {su} 2} + e^{- \frac {su} 2})^2}\right ).$$
Using the substitution $w:=2\cdot \sinh u/2=e^{u/2}-e^{-u/2}=u+\ldots $, one finds that
$$A_1=\mathrm {Res}_{w=0}\frac 1 {w^{2k}} \cdot \frac {\sinh su/2}{(2\cosh su/2 )^2}\cdot \frac  {1}{\cosh u/2}.$$
To show that the polynomial $A_1\in \Q[s]$ is nonzero, we will compute its value for $s=2$ with the help of the addition theorems for $\sinh $ and $\cosh $:
$$\mathrm {Res}_{w=0}\frac 1 {w^{2k}} \cdot \frac {\sinh 2u/2}{(2\cosh 2u/2 )^2}\cdot \frac  {1}{\cosh u/2}$$
$$=\mathrm {Res}_{w=0}\frac 1 {w^{2k}} \cdot  \frac {2\sinh u/2 \cdot \cosh u/2}{4(1 + 2\sinh ^2 u/2)^2}\cdot \frac 1 {\cosh u/2}$$
$$=\mathrm {Res}_{w=0}\frac 1 {w^{2k}} \cdot  \frac {w}{4(1+w^2/2)^2}$$
$$=\text{ coeff. of } w^{2k-2} \text{ in }\frac  {1/4} {(1+w^2/2)^2}\neq 0.$$
Hence, $A_1$ is a nonzero polynomial in $s$. This shows that $A_1$ does not vanish for almost all even integers $s$. It follows that $a(B_c)(\tau)=\pm A$ is a nonzero polynomial in $t$ of degree $1$ for almost all even integers $s$.\proofend 

\section{Proof of the main theorem}\label{section proof}
The proof of the main theorem follows from the previous results by an argument similar to the one in \cite{DG21}. We will focus on the statement on $\mathcal{M}_{sec\geq 0}$; the statement for $\mathcal{M}_{Ric>0}$ is analogous and easier. The main steps are the following (see \cite[Section 6]{DG21} for more details).

As before we will assume that $c$ is odd, $k\geq 2$, $s$ and $t$ are nonzero coprime integers, and $s$ is even. In the following we will fix $c$ and $k$ and choose $s>0$ such that the local datum $a(B_c)(\tau)$ is a nonzero polynomial in $t$ of degree $1$. By Proposition \ref{proposition local datum} there are infinitely many choices for such $s$ and, by Lemma \ref{M cohomology lemma}, different choices for $s$ lead to different manifolds $M_{s,t,c}$, which can be distinguished by their integral cohomology. We will fix a choice for $s$.

By Theorem \ref{theorem finiteness}, the family ${\cal F}_{c,s}=\{ M_{s,t,c}\, \mid \, t \text{ coprime to } s\}$ of $(4k+1)$-dimensional oriented manifolds belongs to finitely many oriented diffeomorphism types. Let us choose a sequence $t_0 <t_1 < t_2 <\ldots$ such that each $M_{s,t_l,c}$ for $l\geq 0$ is diffeomorphic to $M_{s,t_0,c}$ as an oriented manifold and such that the relative $\eta$-invariants $\widetilde{\eta}_\alpha (M_{s,t_l,c})$ for $l\in \N $ are pairwise distinct (see Propositions \ref{proposition eta and local datum} and \ref{proposition local datum}). An orientation-preserving diffeomorphism $M_{s,t_0,c}\to M_{s,t_l,c}$ may not preserve the topological $Spin^c$-structures. However, since $M_{s,t_l,c}$ has only finitely many (namely two) topological $Spin^c$-structures with trivial first Chern class, we may assume, after passing to a subsequence, again denoted by $M_{s,t_l,c}$, that all manifolds in this sequence are diffeomorphic by diffeomorphisms preserving the topological $Spin^c$-structures. Let $M:=M_{s,t_0,c}$ and let $F_l:M\to M_{s,t_l,c}$ be such a diffeomorphism.

Let $g_l := F_l^*(g_{s,t_l,c})$, where $g_{s,t_l,c}$ is the submersion metric of nonnegative sectional and positive Ricci curvature on $M_{s,t,c}$ from Section \ref{section nonnegative curvature}. Since $\eta$-invariants are preserved under pullback, we conclude that the relative $\eta$-invariants of the $Spin^c$-manifold $M$ \wrt \ $g_l$ for $l\in\N $ are pairwise distinct.

Let $\mathcal D $ denote the subgroup of diffeomorphisms of $M$ which preserve its topological $Spin^c$-structure. Note that $\mathcal D $ has finite index in the full diffeomorphism group $\text{Diff}(M)$. Hence, it suffices to show that the elements $[g_{l}]\in \mathcal{R}_{sec\geq0}(M)/ \mathcal D$ for $l\in \N $ defined by $g_l$ represent infinitely many path components.

We argue by contradiction. Suppose there is a path $\tilde{\gamma}:[0,1]\to\mathcal{R}_{sec\geq0}(M)/ \mathcal D$ connecting $[g_{l}]$ to $[g_{l^\prime}]$ with $l\neq l^\prime $. By Ebin's slice theorem \cite{Eb70}, this path can be lifted to a continuous path $\gamma$ in $\mathcal{R}_{sec\geq 0}(M)$ with $\gamma(0)=g_{l}$ and $\gamma(1)=\Phi^*(g_{l^\prime })$ for some $\Phi\in\mathcal D$. Since $\eta$-invariants are preserved under pullback, it follows that the relative $\eta$-invariants of the $Spin^c$-manifold $M$ \wrt \ $\gamma(0)=g_{l}$ and $\gamma(1)=\Phi^*(g_{l^\prime})$ are distinct.

The path $\gamma $ may be deformed inside of $\mathcal{R}_{scal >0}(M)$ to a path $\hat{\gamma}$ with the same endpoints as $\gamma$ and whose interior points lie in $\mathcal{R}_{Ric>0}(M)$ (this can be done via Ricci flow using \cite{BW07}). Since the relative $\eta$-invariant is constant on path components of $\mathcal{R}_{scal>0}(M)$ (see \cite[p. 417]{APSII75}, \cite[Proposition 3.3]{DG21}), we get a contradiction.

Hence, the classes $[g_{l}]$ for $l\in\N $ represent infinitely many pairwise distinct path components of $\mathcal{R}_{sec\geq0}(M)/ \mathcal D$. Since $\mathcal D $ has finite index in $\text{Diff}(M)$, the same holds for the moduli space $\mathcal{M}_{sec\geq 0}(M)$. As explained in the beginning, we can argue in this way for infinitely many choices of $s$. Hence, we obtain infinitely many manifolds $M_i:=M_{s_i,t_i,c}$, indexed by $i\in \N $, which can be distinguished by their integral cohomology, such that for each $i\in \N $ the moduli space $\mathcal{M}_{sec\geq 0}(M_i)$ has infinitely many path components. This completes the proof of the first statement of the main theorem. An analogous argument gives the statement for $\mathcal{M}_{Ric>0}$.\proofend

\begin{remark}\label{dim five remark}  For $k=1$, the manifolds $\overline M_{s,t,c}$ given in Definition \ref{simply connected definition} can be shown to be diffeomorphic to $S^2\times S^3$ (for $c$ odd, $s$ even, and $s$ and $t$ coprime). The $\Z _2$-quotients $M_{s,t,c}$ can be described as total spaces of $S^1$-principal bundles over $B_c\cong \C P^2 \sharp -\C P^2 $ and fall into finitely many diffeomorphism types. Their moduli spaces of metrics of nonnegative sectional curvature and positive Ricci curvature also have infinitely many path components (see also the recent work of Goodman and Wermelinger \cite{GW22,W20} on such moduli spaces for the class of all orientable nonspin $\Z _2$-quotients of $S^2\times S^3$).
\end{remark}

\end{document}